\documentclass[12pt, a4paper, twoside]{article}
\usepackage{indentfirst}
\usepackage{ifpdf}
\usepackage{graphicx,natbib,amssymb,lineno}
\small
\bibliography{mybibfile}
\ifpdf
\usepackage[%
  pdftitle={ title{The Art of Space Filling in Penrose Tilings and Fractals},%
  pdfauthor={San Le},%
  pdfkeywords={instructions for use, elsart, document class},%
  pdfstartview=FitH,%
  bookmarks=true,%
  bookmarksopen=true,%
  breaklinks=true,%
  colorlinks=true,%
  linkcolor=blue,anchorcolor=blue,%
  citecolor=blue,filecolor=blue,%
  menucolor=blue,pagecolor=blue,%
  urlcolor=blue]{hyperref}
\else
\usepackage[%
  breaklinks=true,%
  colorlinks=true,%
  linkcolor=blue,anchorcolor=blue,%
  citecolor=blue,filecolor=blue,%
  menucolor=blue,pagecolor=blue,%
  urlcolor=blue]{hyperref}
\fi

\makeatletter
\def\elsartstyle{%
    \def\normalsize{\@setfontsize\normalsize\@xiipt{14.5}}
    \def\small{\@setfontsize\small\@xipt{13.6}}
    \let\footnotesize=\small
    \def\large{\@setfontsize\large\@xivpt{18}}
    \def\Large{\@setfontsize\Large\@xviipt{22}}
    \skip\@mpfootins = 18\p@ \@plus 2\p@
    \normalsize
}
\@ifundefined{square}{}{}
\makeatother

\NeedsTeXFormat{LaTeX2e}

\usepackage{afterpage}
\usepackage{amssymb}
\usepackage{amsmath}
\usepackage{subeqnarray}
\usepackage{setspace}
\newsavebox{\astrutbox}
\sbox{\astrutbox}{\rule[-5pt]{0pt}{20pt}}

\pagestyle{plain}
\begin{document}

\title{The Art of Space Filling in Penrose Tilings and Fractals}


\medskip
\begin{center}
{\Large\bfseries The Art of Space Filling in Penrose Tilings and Fractals} \par\medskip
{ By SAN LE}\par
\medskip
{\small California, USA\\[3pt]}\par
\medskip
\end{center}
\medskip


\begin{abstract}
\raggedright 
{
%
%
%
%
\parindent=20pt 
%
%
%
%
}

Incorporating designs into the tiles that form tessellations presents an interesting
challenge for artists.  Creating a viable MC Escher like image that works esthetically as well
as functionally requires resolving incongruencies at a tile's edge while constrained
by its shape.  Escher was the most well known practitioner in this style
of mathematical visualization, but there are significant mathematical shapes
to which he never applied his artistry.  These shapes can incorporate designs that form
images as appealing as those produced by Escher, and our paper explores this for
traditional tessellations, Penrose Tilings, fractals, and fractal/tessellation
combinations.  To illustrate the versatility of tiling art, images were created
with multiple figures and negative space leading to patterns distinct from
the work of others.  
\footnote[1]
{ 
\raggedright 
Received June, 2011

2010 Mathematics Subject Classification. Primary 00A66, 05B45, 97M80; Secondary 52C20, 28A80.

\emph{Key words and phrases.} Space Filling; Penrose Tilings; Escher; Symmetry; Tilings; Fractals; Tessellation

}

\end{abstract}



\section{Introduction}
\raggedright 
\parindent=20pt 

MC Escher was the most prominent artist working with tessellations and space filling.
Forty years after his death, his creations are still foremost in people's minds in the
field of tiling art.  One of the reasons Escher continues to hold such a monopoly
in this specialty are the unique challenges that come with creating Escher type designs
inside a tessellation\cite[]{Raed2005}.  When an image is drawn into a tile
and extends to the tile's edge, it
introduces incongruencies which are resolved by continuously aligning and
refining the image.  This is particularly true when the image consists of the
lizards, fish, angels, etc. which populated Escher's tilings because they do not have
the quadrilateral symmetry that would make it possible to arbitrarily rotate the image $\pm$ 90, 180
degrees and have all the pieces fit\cite[]{James96}.  Rather, they have bilateral symmetry which requires creating
a compliment for every edge.  This is true for any type of tile that incorporates
such an image.

A collection of papers in honor of Escher, \emph{MC Escher's Legacy:
A Centennial Celebration} contains a comprehensive study of his work \cite[]{Schat05}.
Of the articles emphasizing art in this collection, most authors produced tile images
that continued the practice of having one dominant figure in a tile that is
completely filled.  There are greater possibilities beyond this, and this paper
contains the result of our studies.  

The rules to creating tiling art are straightforward, and we describe the process and apply it to
the mathematical geometries that Escher did not have a chance to include in his vast
body of work.  In particular, we create tiling art based on constructions that have
gained prominence since Escher's time: Penrose Tilings and fractals.  In the case of
the latter, we also put a tessellation inside a fractal tile to allow for growth in various
directions.  To give our designs a classical quality, they will consist of intertwined
human figures.  We also deviate from Escher and others by having
negative spaces between the figures that allowed for different combinations in connecting the tiles.

Most of the visualizations related to Penrose Tilings and fractals are abstract representations
purely in service of the underlying numbers \cite[]{Bandt97}.  This represents a missed opportunity as
the geometry of tessellations and fractals is highly suited to the application of
design art offering a rich and mostly unmined domain awaiting exploration.  There
are rewards in its provinces beyond what can be achieved otherwise.  The images produced will often
excite that part of the mind which responds to the inherent esthetic of mathematical
shapes \cite[]{Rota77, Abas76}.  In addition, the artist enjoys another
level of discovery to his or her creation.
Unlike a picture in isolation, the tiling images in this paper form patterns that only revealed themselves
once all the tiles are assembled.  The compelling nature of this art was best
explained by Escher himself, saying, ".. For once one
has crossed over the threshold of the early stages this activity
takes on a more worth than any other form of decorative art" \cite[]{Ern95}.

\section{Creating a tessellation}
\raggedright 
\parindent=20pt

\begin{figure}
  \centerline{\includegraphics[height=3.0in,width=3.0in,angle=0]{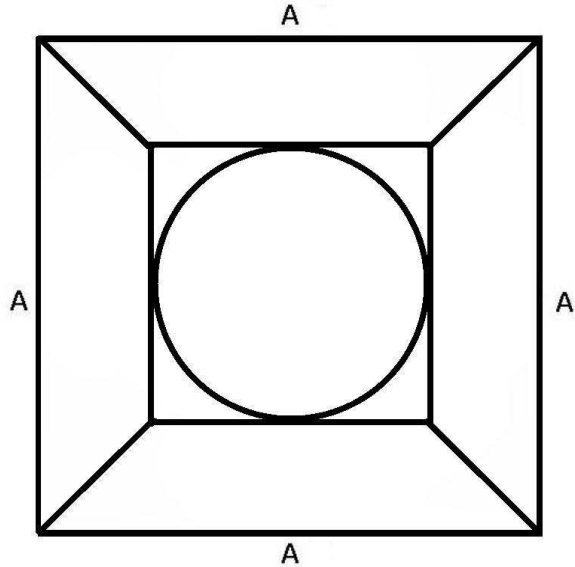}}
  \caption{Simple tile with quadrilateral symmetry}
  \label{square_quad_symm2}
  \end{figure}

\begin{figure}
  \centerline{\includegraphics[height=4.94in,width=5.0in,angle=0]{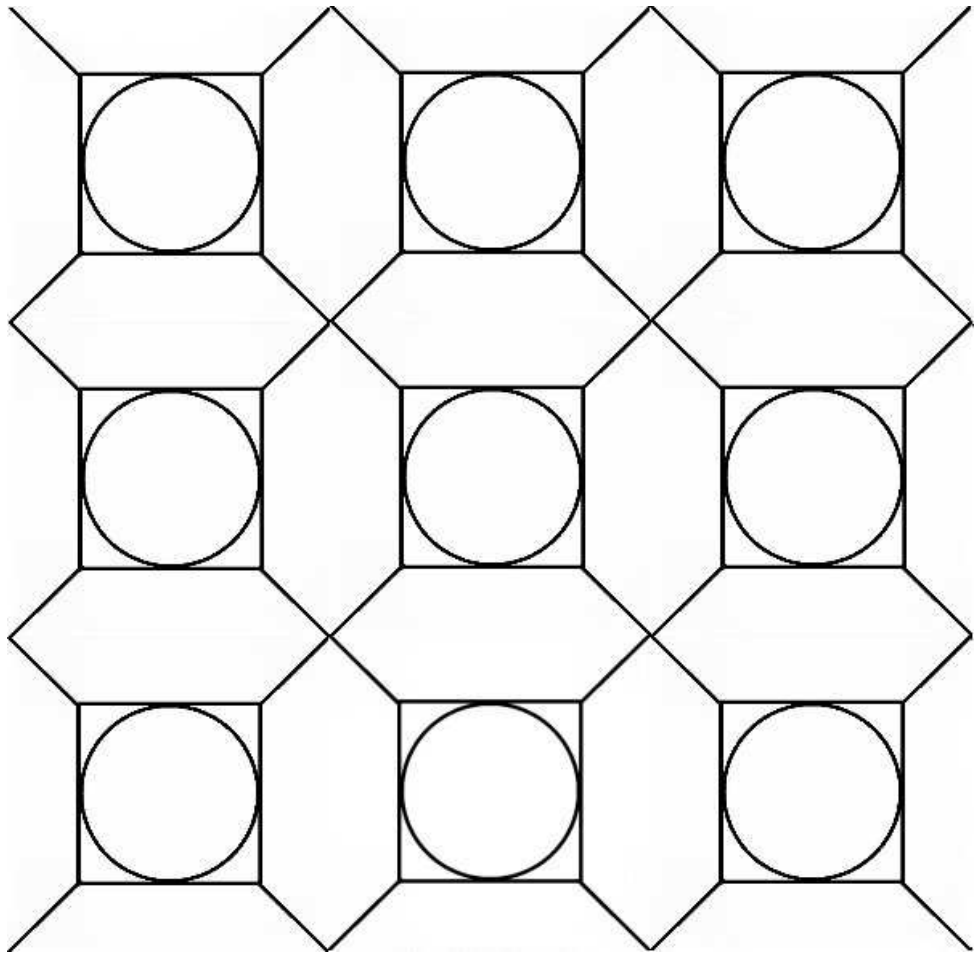}}
  \caption{Tessellation of Figure \ref{square_quad_symm2} }
  \label{square_quad_symm_tessel}
  \end{figure}

To illustrate the process and challenges of making tiles that contain bilateral
symmetry, we begin with the simplest tessellation made of squares.  Figure \ref{square_quad_symm2}
illustrates a square with quadrilateral symmetry, and any $\pm$ 90, 180 degree rotation will
give back the same image.  Any side \emph{A} can be connected to another side \emph{A}
and the resulting tessellation is given in Figure \ref{square_quad_symm_tessel}. 

%
%
%
%
%

\begin{figure}
  \centerline{\includegraphics[height=3.0in,width=3.0in,angle=0]{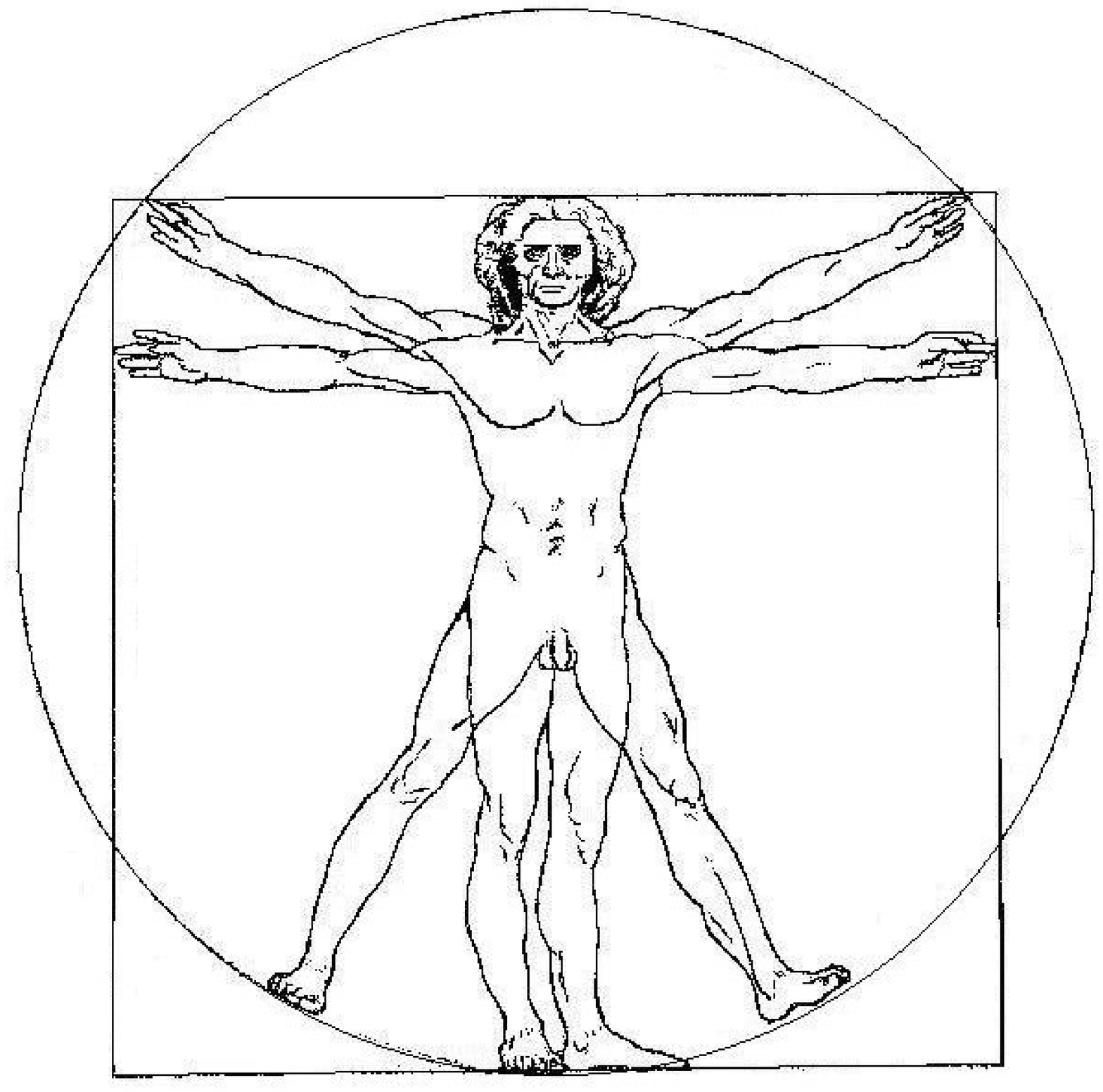}}
  \caption{Leonardo da Vinci's \emph{The Vitruvian Man}}
  \label{Da_Vinci_Vitruve_Luc_Viatour_small2f_crop}
  \end{figure}

\begin{figure}
  \centerline{\includegraphics[height=2.0in,width=2.0in,angle=0]{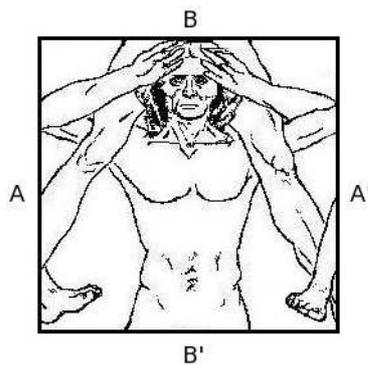}}
  \caption{A single tile of \emph{The Vitruvian Man} for tessellation}
  \label{Da_Vinci_Vitruve_Luc_Viatour_cell_border2}
  \end{figure}

\begin{figure}
  \centerline{\includegraphics[height=5.44in,width=4.0in,angle=0]{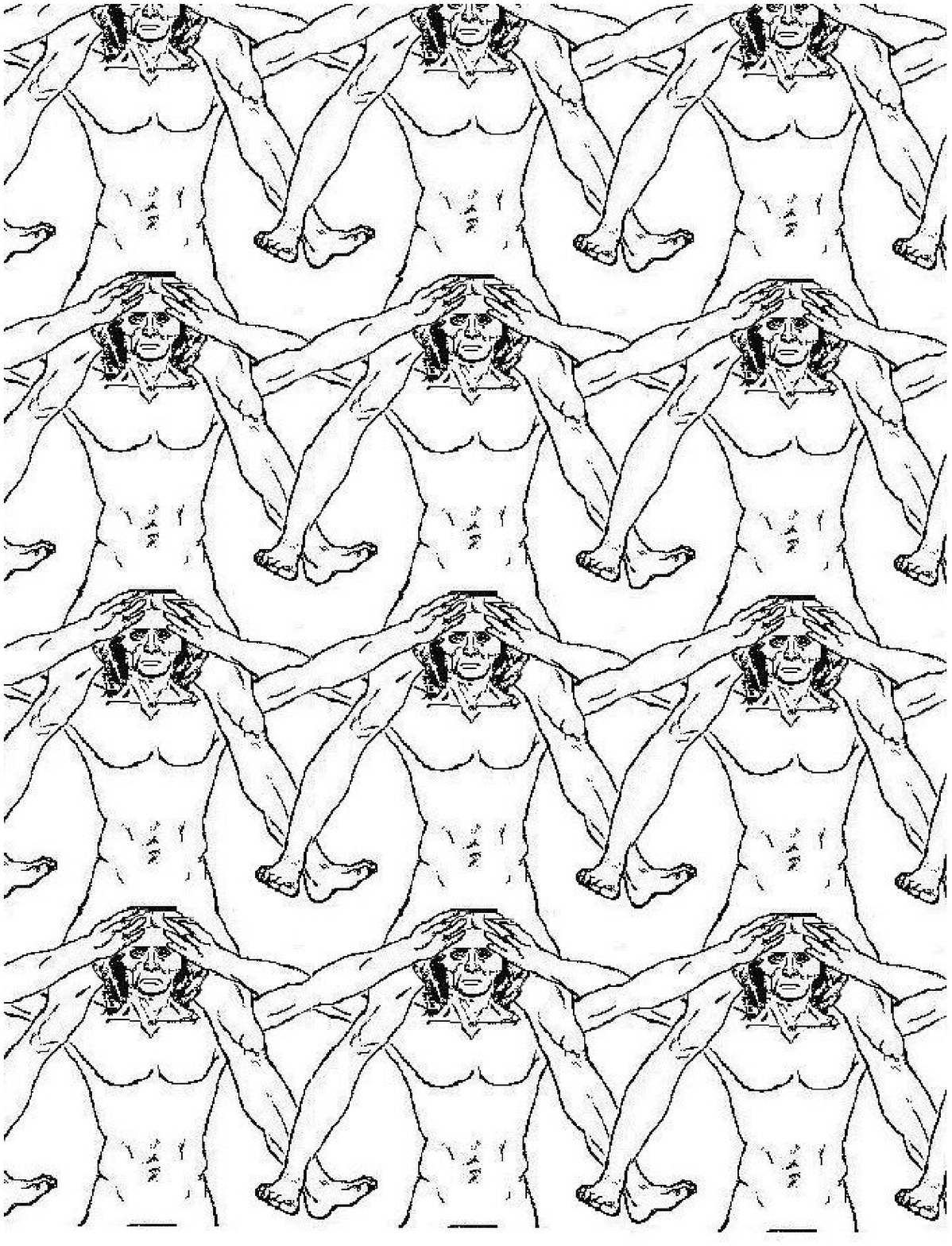}}
  \caption{Tessellation of \emph{The Vitruvian Man} }
  \label{Da_Vinci_Vitruve_Luc_Viatour_tessel}
  \end{figure}

In contrast, when the design is of a person as in \emph{Leonardo da Vinci}'s \emph{The Vitruvian Man}
shown in Figure \ref{Da_Vinci_Vitruve_Luc_Viatour_small2f_crop}\cite[]{Mir95}, it is not possible
to match every side to every other if the image is extended to the edges.  The
overlapping of arms and legs requires that each side is complimented by another side.
To create a tessellation from Figure \ref{Da_Vinci_Vitruve_Luc_Viatour_small2f_crop},
we label sides as before, but now we have four unique sides.
In the single tile given by Figure \ref{Da_Vinci_Vitruve_Luc_Viatour_cell_border2},
side \emph{A} connects to side \emph{A'} and \emph{B} to \emph{B'}.  The resulting tessellation is
given in Figure \ref{Da_Vinci_Vitruve_Luc_Viatour_tessel}.

\section{Beyond Escher}
\raggedright 
\parindent=20pt 

Once the simple rules of making a tessellation are understood, other possibilities exist
in how tiles are connected and the types of images they contain.

\subsection{Other Images and Connecting Rules}
\raggedright 
\parindent=20pt 

\begin{figure}
  \centerline{\includegraphics[height=7.96in,width=6.0in,angle=0]{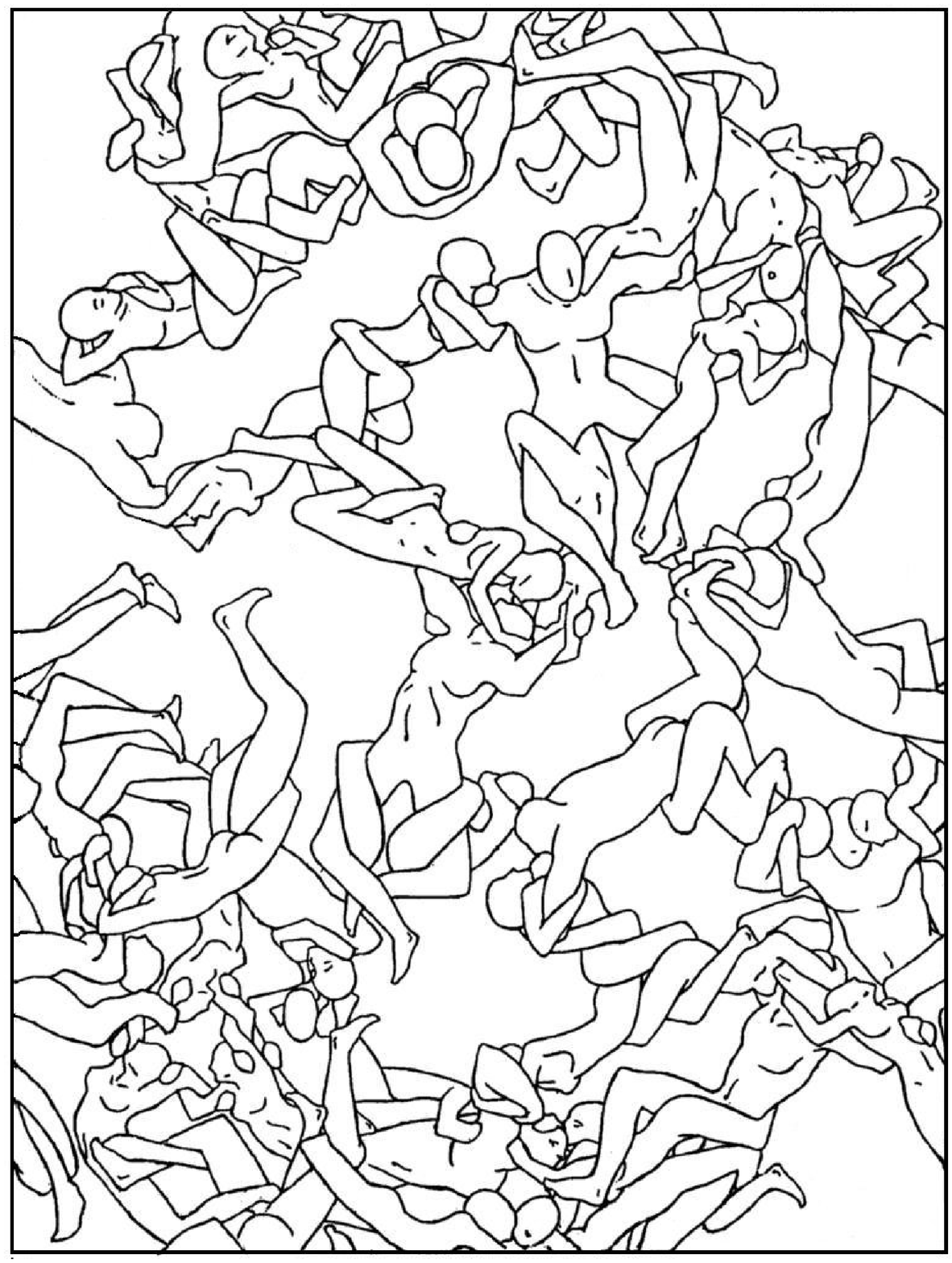}}
  \caption{A tile made up of multiple figures}
  \label{web_cell_border}
  \end{figure}

\begin{figure}
  \centerline{\includegraphics[height=7.86in,width=6.0in,angle=0]{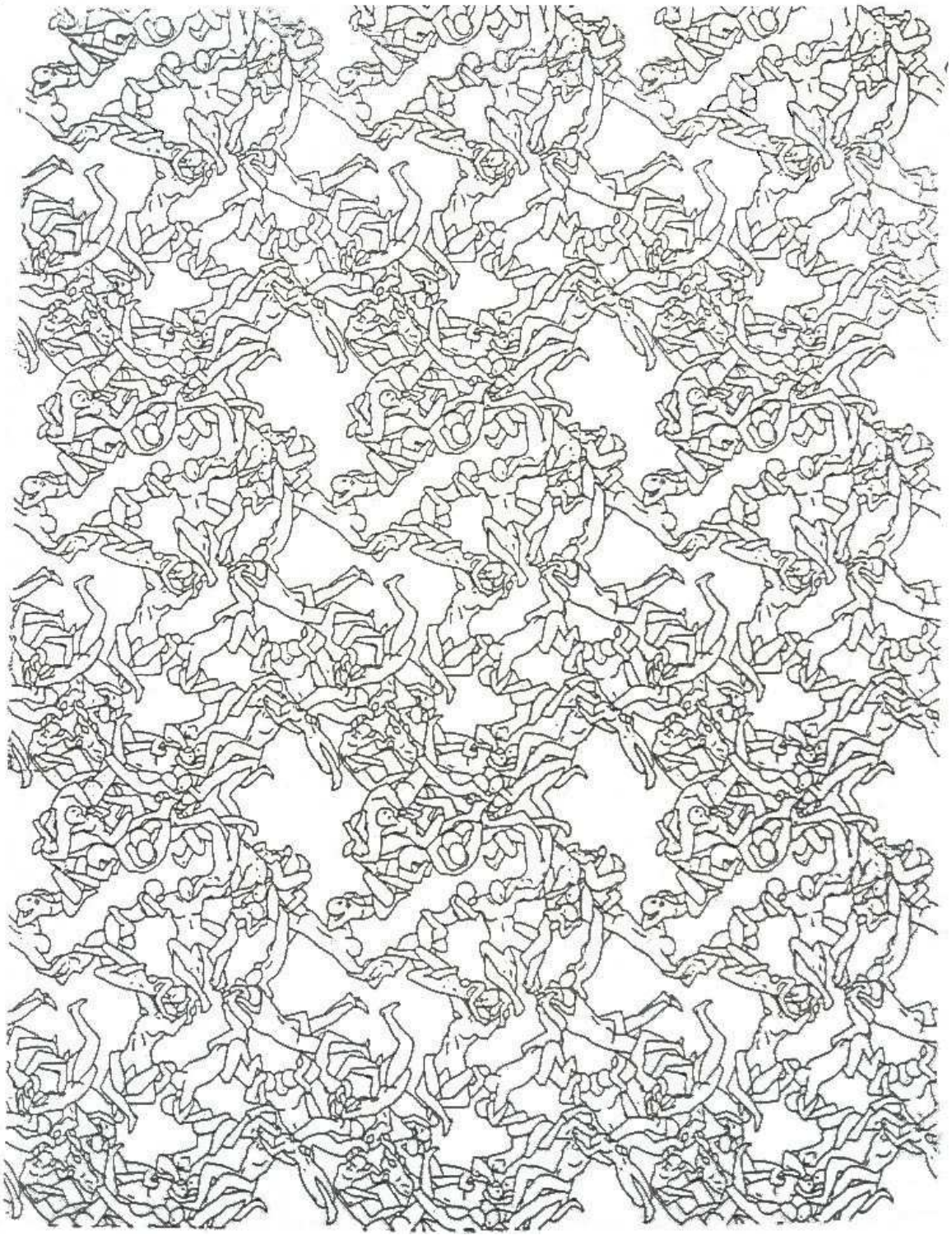}}
  \caption{Tessellation of Figure \ref{web_cell_border} }
  \label{web_clean}
  \end{figure}

Instead of limiting a tile to one prominent entity as typically seen in most work emulating
Escher, Figure \ref{web_cell_border} contains multiple subjects and its rectangular
tile leads to the tessellation of Figure \ref{web_clean}.

\begin{figure}
  \centerline{\includegraphics[height=3.0in,width=3.0in,angle=0]{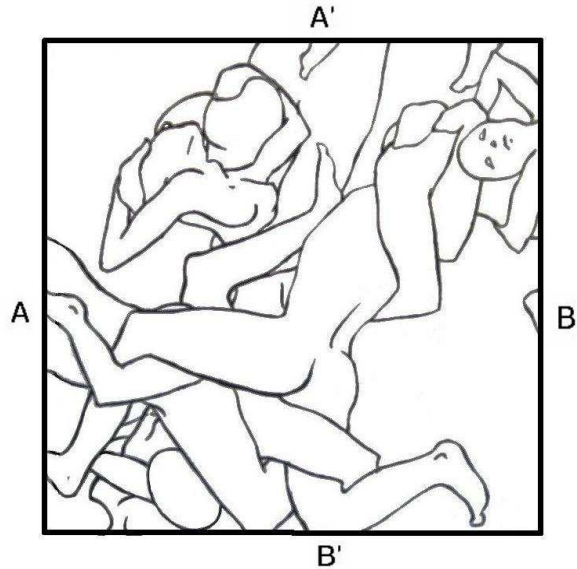}}
  \caption{ Tile having adjacent complimentary sides }
  \label{circle_cellB_borderB}
  \end{figure}

\begin{figure}
  \centerline{\includegraphics[height=5.03in,width=5.0in,angle=0]{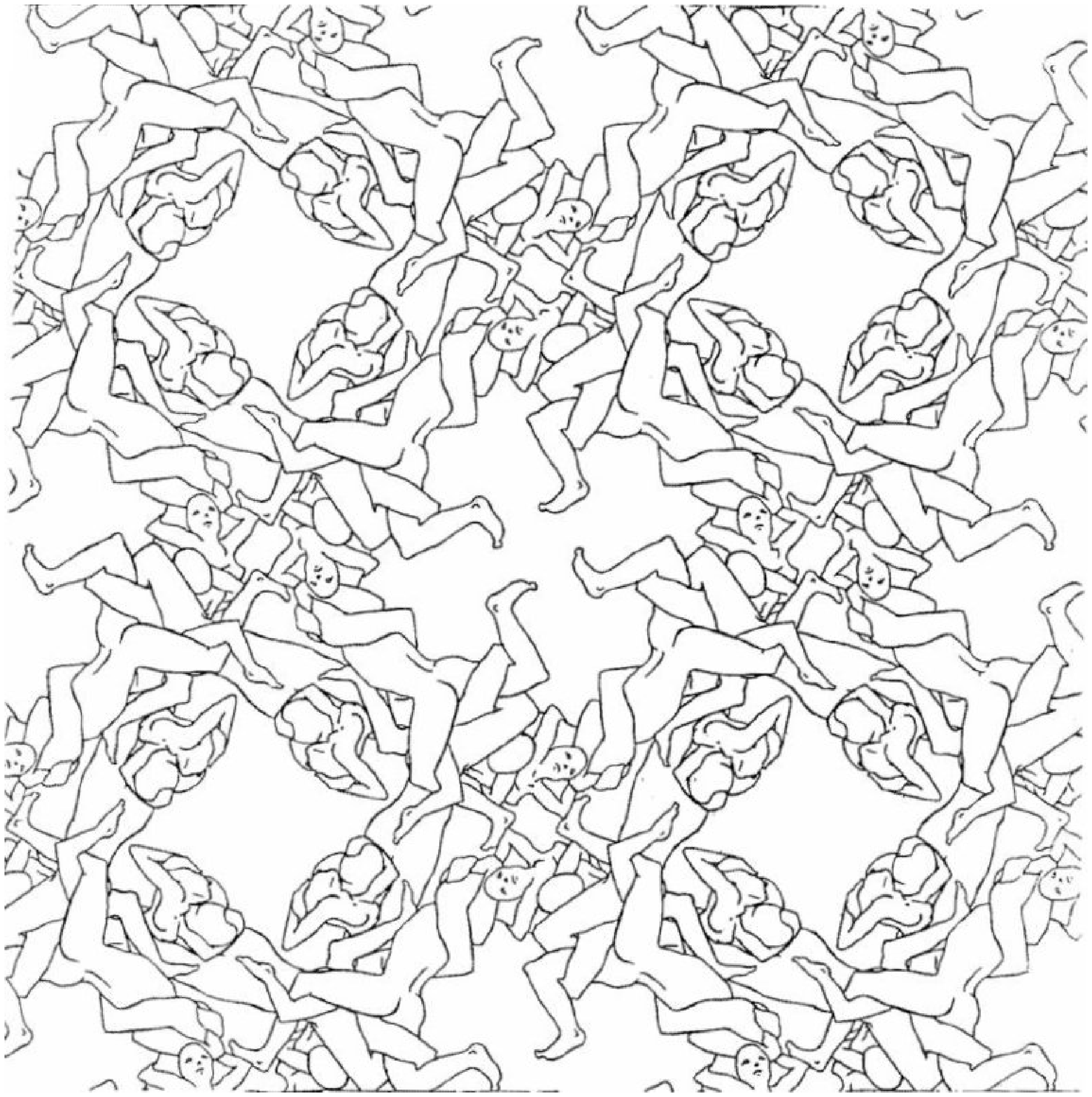}}
  \caption{Tessellation of Figure \ref{circle_cellB_borderB} }
  \label{circle_crop}
  \end{figure}

\begin{figure}
  \centerline{\includegraphics[height=3.0in,width=3.0in,angle=0]{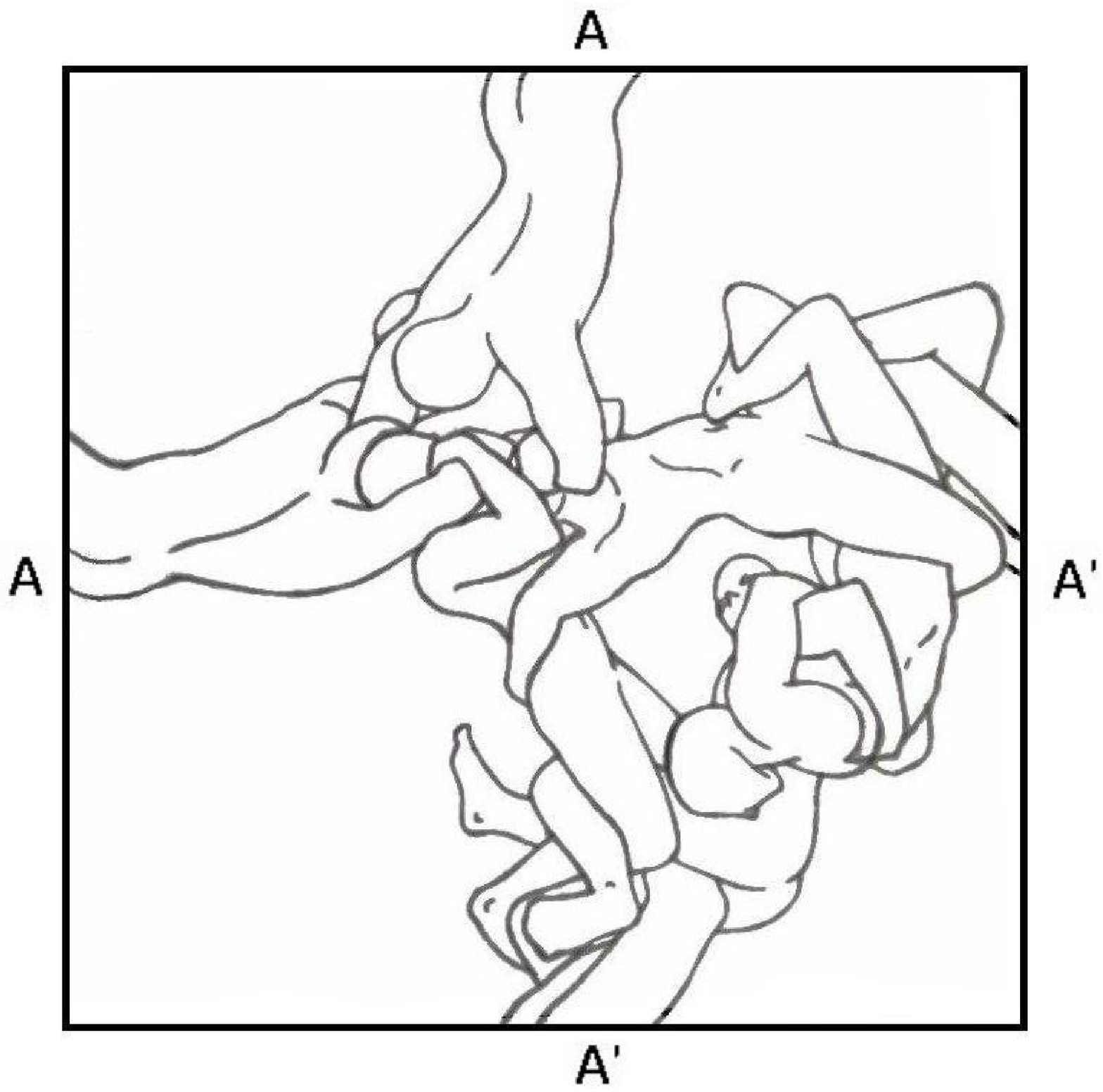}}
  \caption{ Tile having two adjacent complimentary sides }
  \label{square_cellB_borderC}
  \end{figure}

\begin{figure}
  \centerline{\includegraphics[height=5.07in,width=5.0in,angle=0]{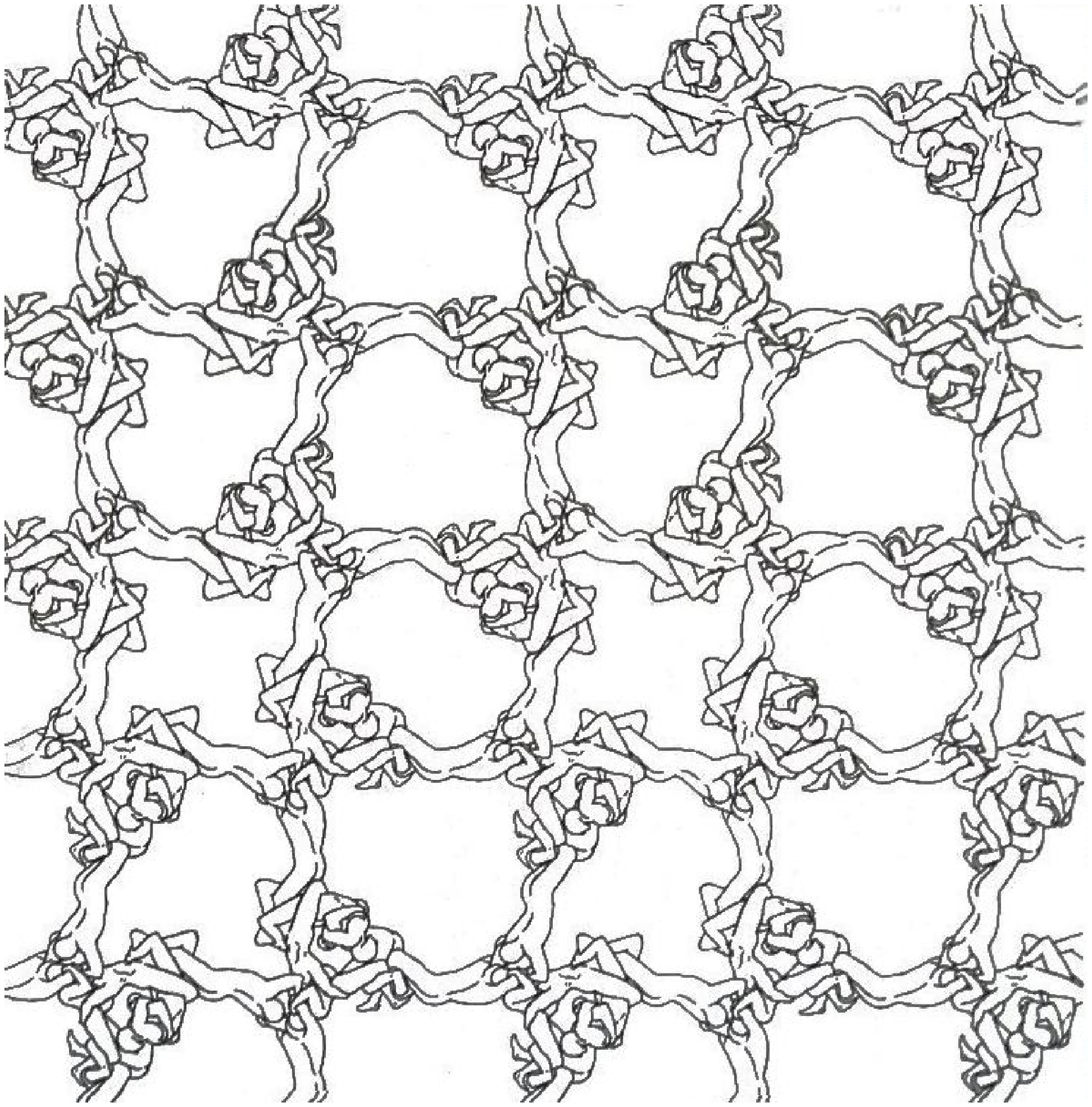}}
  \caption{Tessellation of Figure \ref{square_cellB_borderC} }
  \label{square_crop}
  \end{figure}

The next two figures show other matching rules that square tiles can have.
As seen in Figure \ref{circle_crop} formed by the tiles of Figure \ref{circle_cellB_borderB},
by making the complimentary sides adjacent, the final image forms an interesting pattern of
swirling figures rather than a simple shifting of the tile vertically or horizontally.  In
Figure \ref{square_cellB_borderC}, the same adjacent compliments exist, but with each side
complimenting two sides instead of one.  This leads to the image seen in Figure \ref{square_crop}
which has the same turning of the figures, but also allows for non periodic
arrangements.  Once one row of tiles is set, the next row has two possible arrangements
depending on how the first tile of the new row is placed.  This is true for every subsequent row.


%

%

%

\subsection{Penrose Tilings}
\raggedright 
\parindent=20pt

\begin{figure}
  \centerline{\includegraphics[height=4.65in,width=3.0in,angle=0]{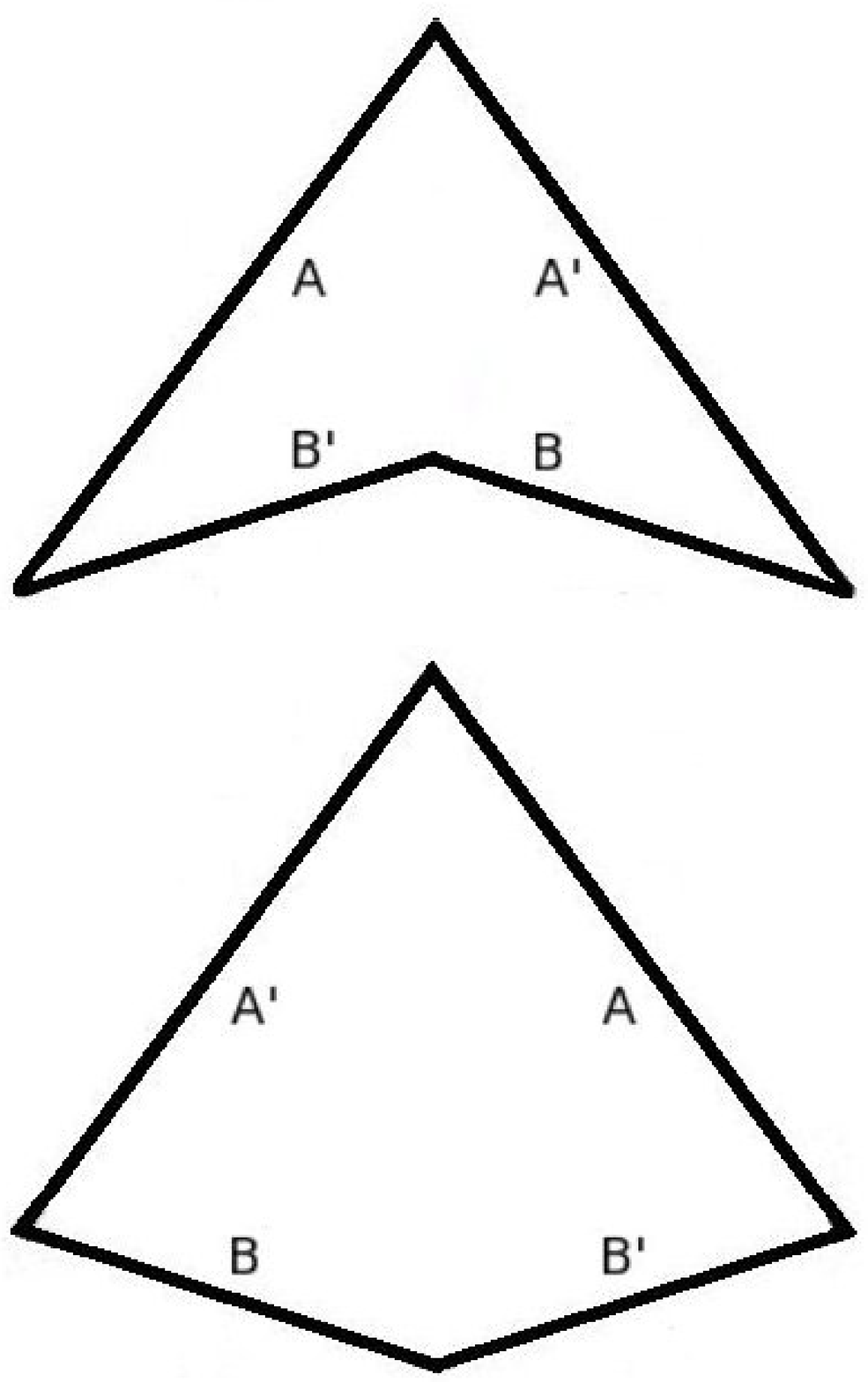}}
  \caption{ Kite and Dart Penrose Tilings }
  \label{Penrose_Dart_label3}
  \end{figure}

\begin{figure}
  \centerline{\includegraphics[height=5.41in,width=5.0in,angle=0]{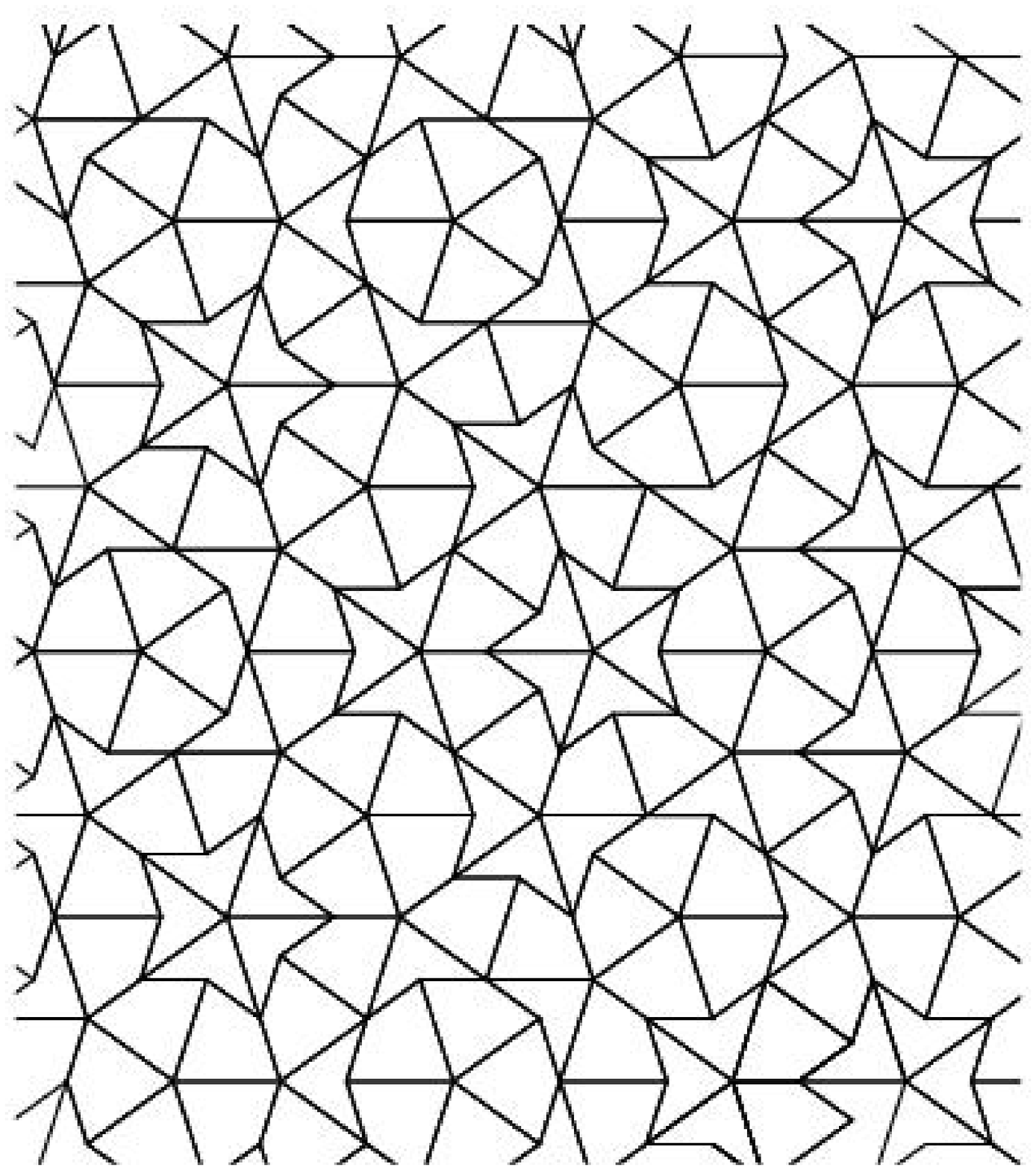}}
  \caption{ Tessellation of Figure \ref{Penrose_Dart_label3} }
  \label{Kite_and_Dart}
  \end{figure}

\begin{figure}
  \centerline{\includegraphics[height=6.76in,width=5in,angle=0]{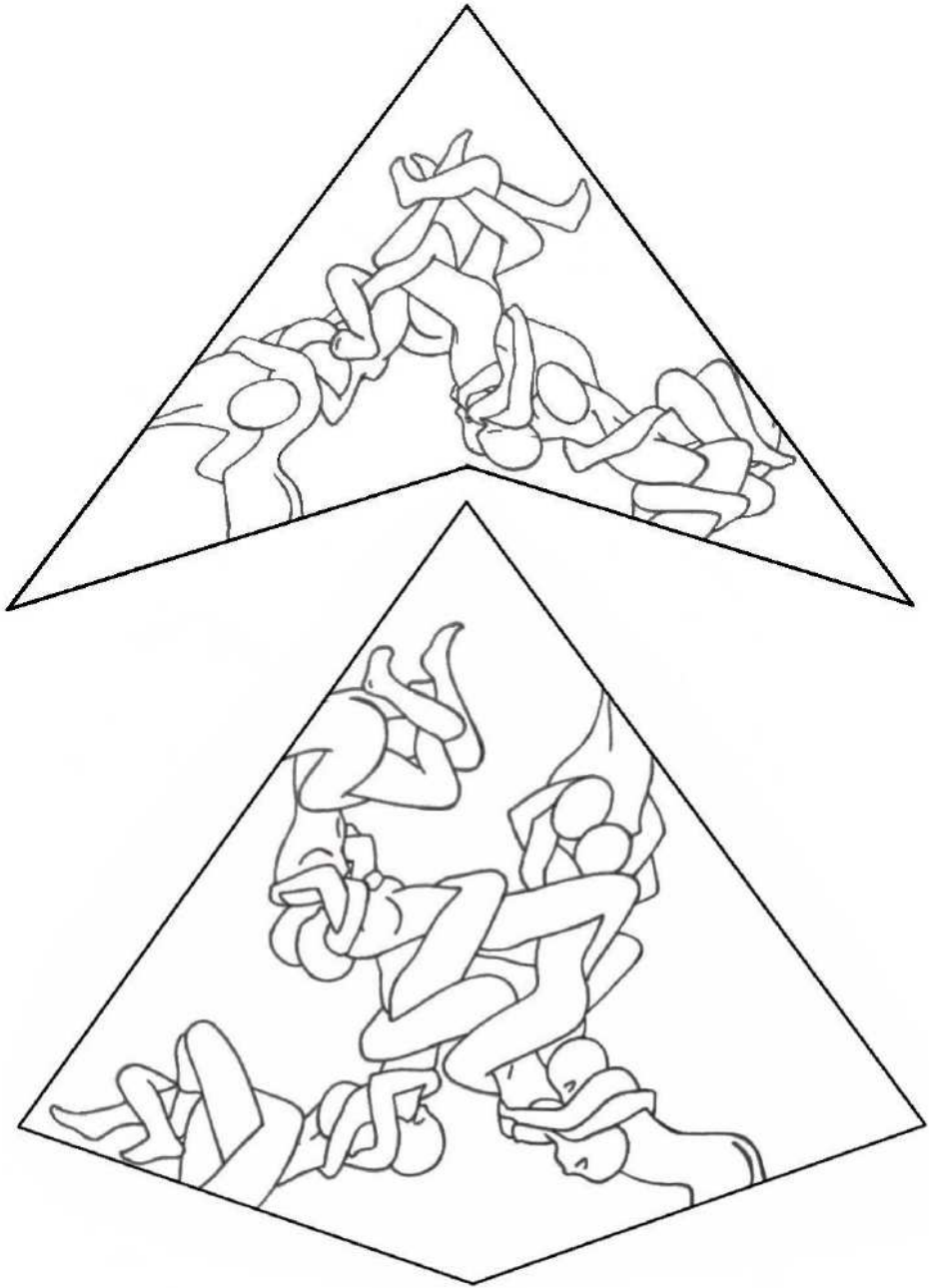}}
  \caption{ Space filled kite and dart Penrose Tiles }
  \label{penrose_Dart_cells}
  \end{figure}

\begin{figure}
  \centerline{\includegraphics[height=8.0in,width=6in,angle=0]{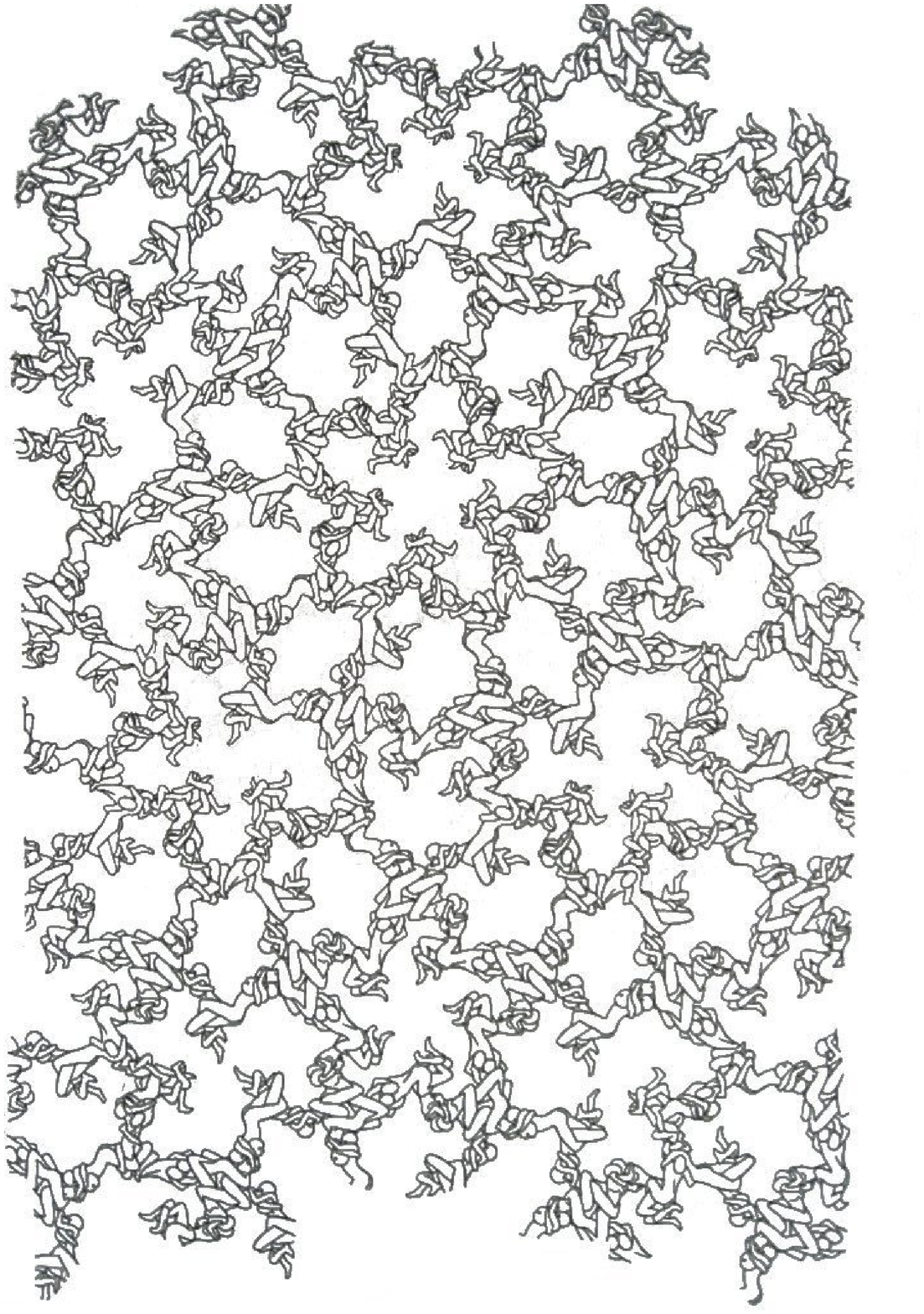}}
  \caption{ Tessellation of Figure \ref{penrose_Dart_cells} }
  \label{penrose2_paper}
  \end{figure}

Though Escher was friends with Roger Penrose, he died before applying his space filling work to Penrose's
aperiodic tilings\cite[]{Dun94, Kaplan04}.  These types of tilings are even more challenging to incorporating designs
because they are made of two tiles which, depending on the arrangement, may not form a
repeating pattern.  Fortunately, the more interesting and non-repeating (hence aperiodic) arrangements
are created by following side matching rules in the same manner as discussed
before\cite[]{Lord91}.  The Kite and Dart tiles, given in Figure \ref{Penrose_Dart_label3}, show how every side
is connected to create the aperiodic tiling of \ref{Kite_and_Dart}.  Using these same rules,
the artistic tiles in Figure \ref{penrose_Dart_cells} lead to \ref{penrose2_paper}.

\begin{figure}
  \centerline{\includegraphics[height=3.33in,width=3.0in,angle=0]{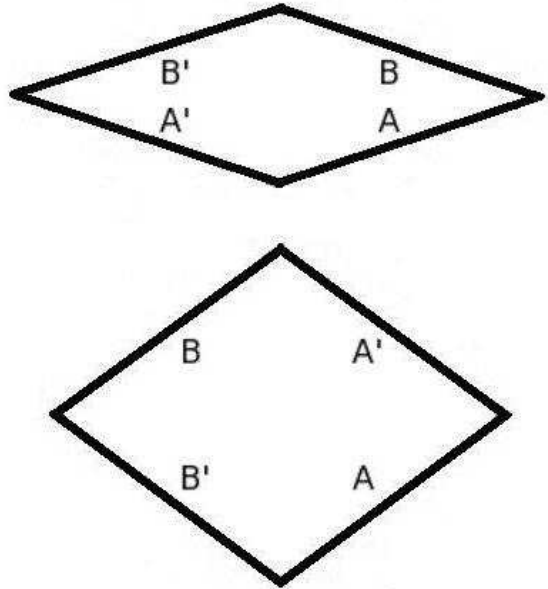}}
  \caption{ Rhombus Penrose Tilings }
  \label{Penrose_Rhombi_label2}
  \end{figure}

\begin{figure}
  \centerline{\includegraphics[height=5.41in,width=5.0in,angle=0]{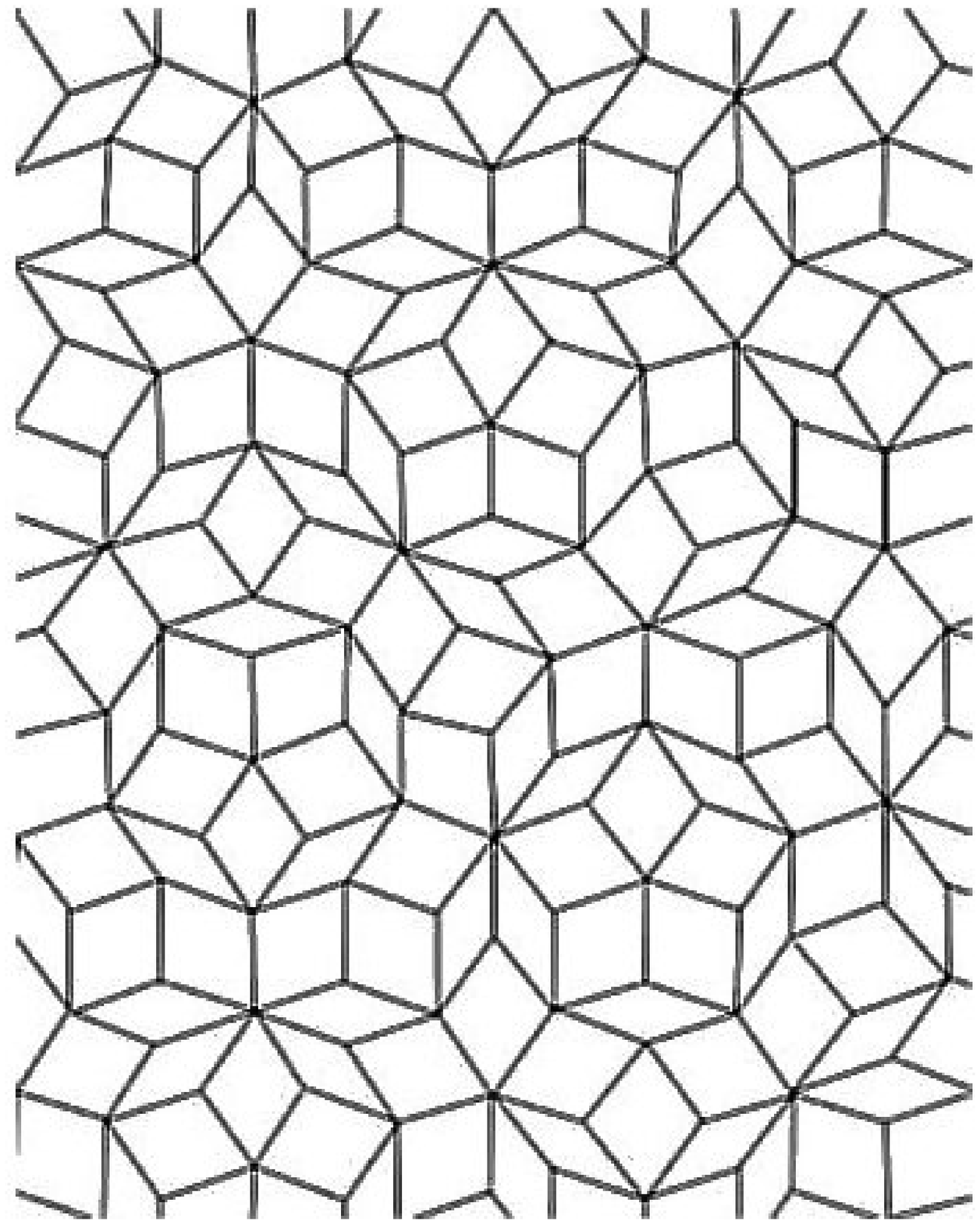}}
  \caption{ Tessellation of Figure \ref{Penrose_Rhombi_label2} }
  \label{penrose_rhombus3c}
  \end{figure}

\begin{figure}
  \centerline{\includegraphics[height=5.23in,width=5in,angle=0]{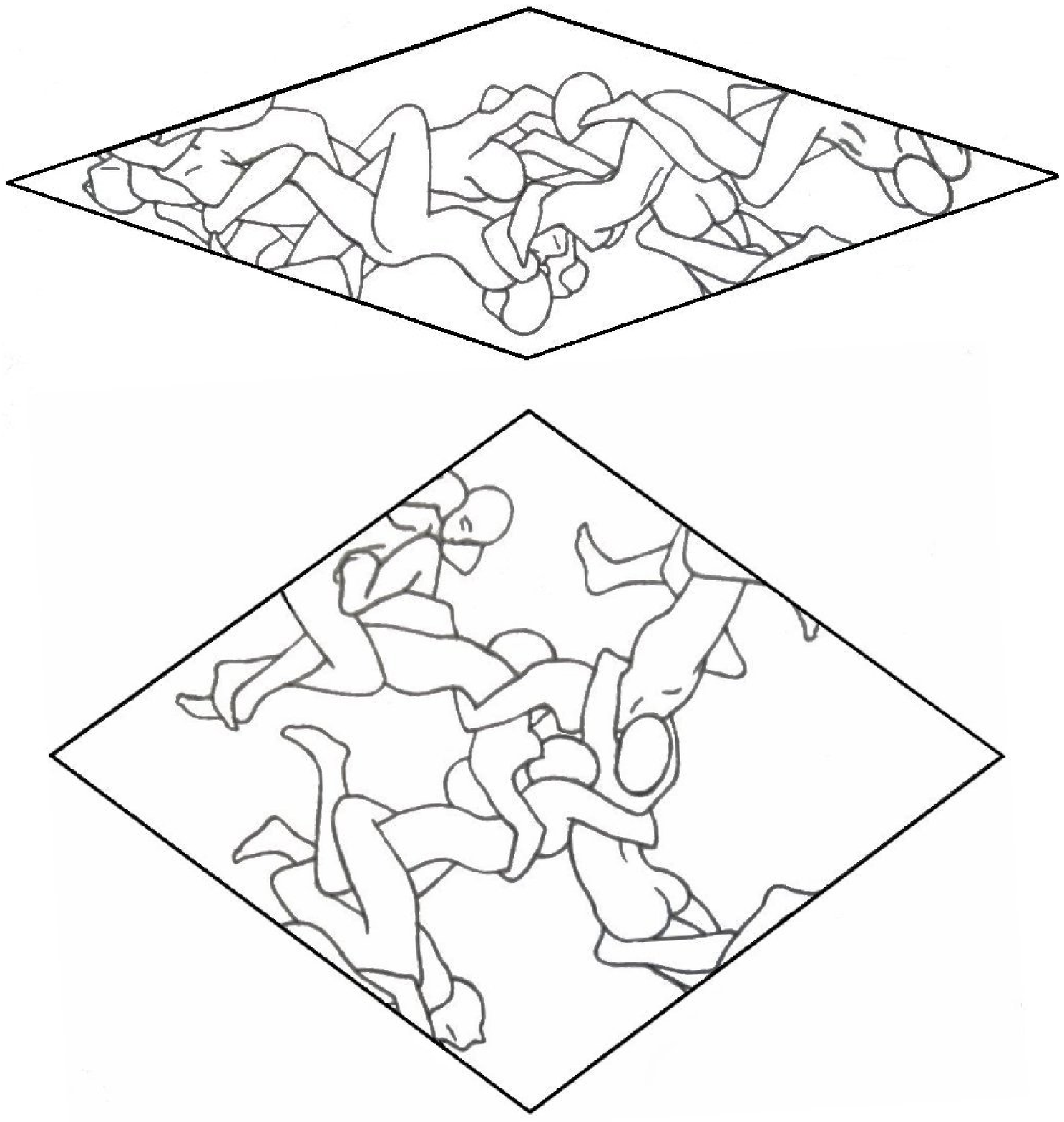}}
  \caption{ Space filled rhombus Penrose Tiles }
  \label{penrose_rhombus_cells}
  \end{figure}

\begin{figure}
  \centerline{\includegraphics[height=8in,width=6in,angle=0]{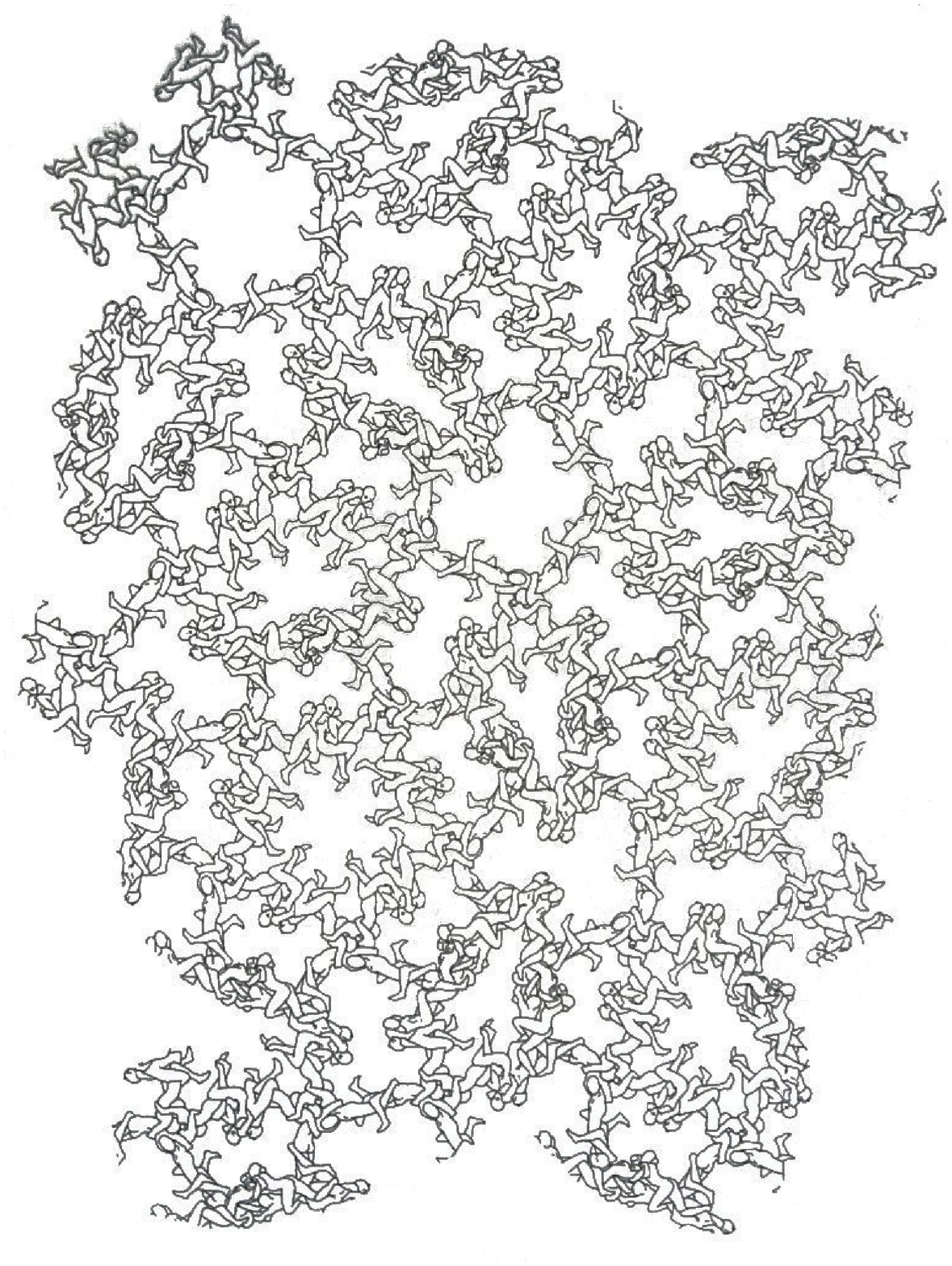}}
  \caption{ Tessellation of Figure \ref{penrose_rhombus_cells} }
  \label{penrose_paper}
  \end{figure}

The other Penrose Tiling made of fat and thin rhombi has similar rules for creating
an aperiodic pattern.  Figure \ref{Penrose_Rhombi_label2} shows how the sides match
to create \ref{penrose_rhombus3c}.  The individual tiles containing a design and
the resulting image are given in Figures \ref{penrose_rhombus_cells} and \ref{penrose_paper}
respectively.

\subsection{Fractal Tilings}
\raggedright 
\parindent=20pt

\begin{figure}
  \centerline{\includegraphics[height=4.32in,width=4in,angle=0]{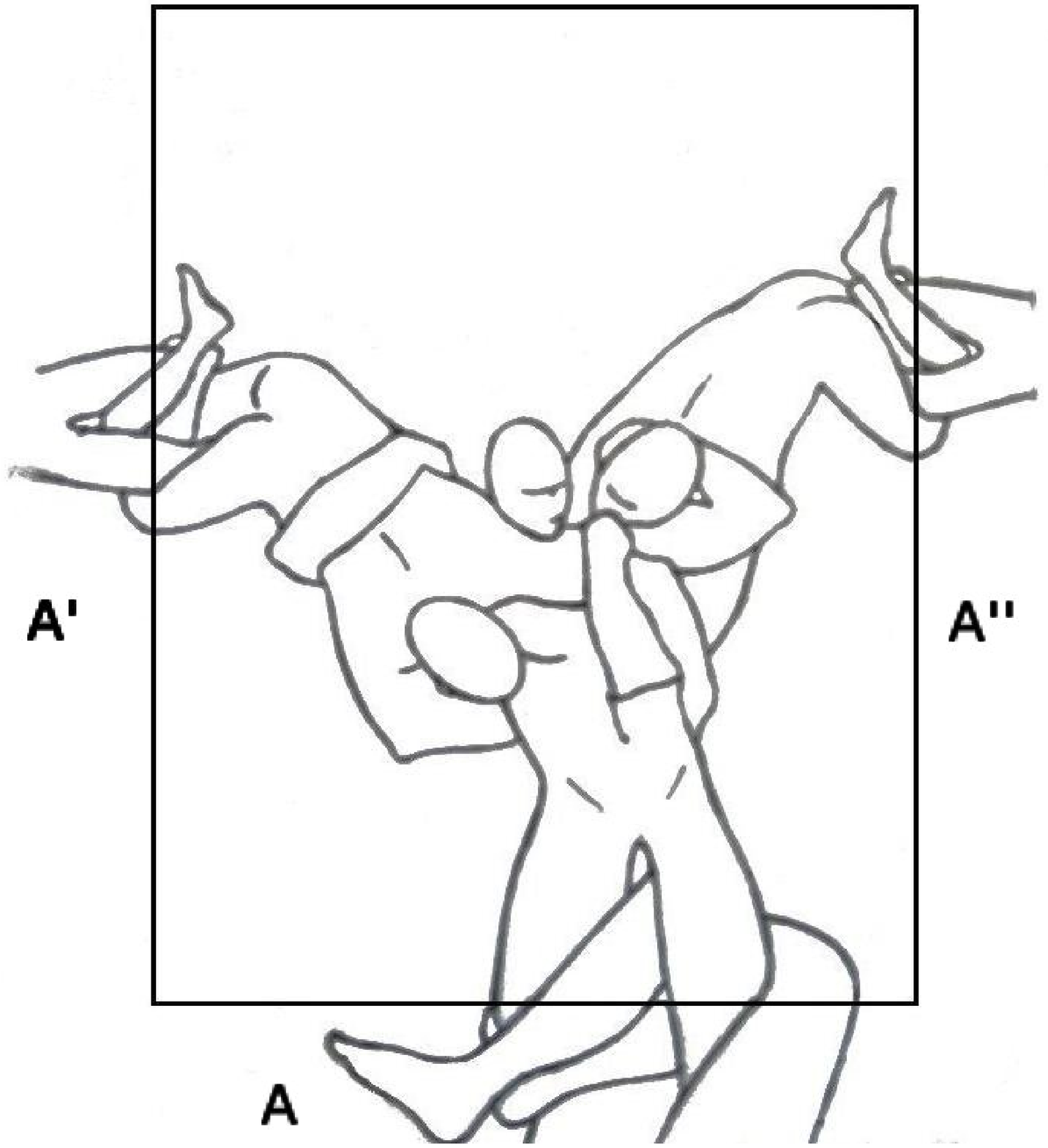}}
  \caption{ Fractal tile }
  \label{cell2D_crop_border2}
  \end{figure}

\begin{figure}
  \centerline{\includegraphics[height=3.33in,width=5in,angle=0]{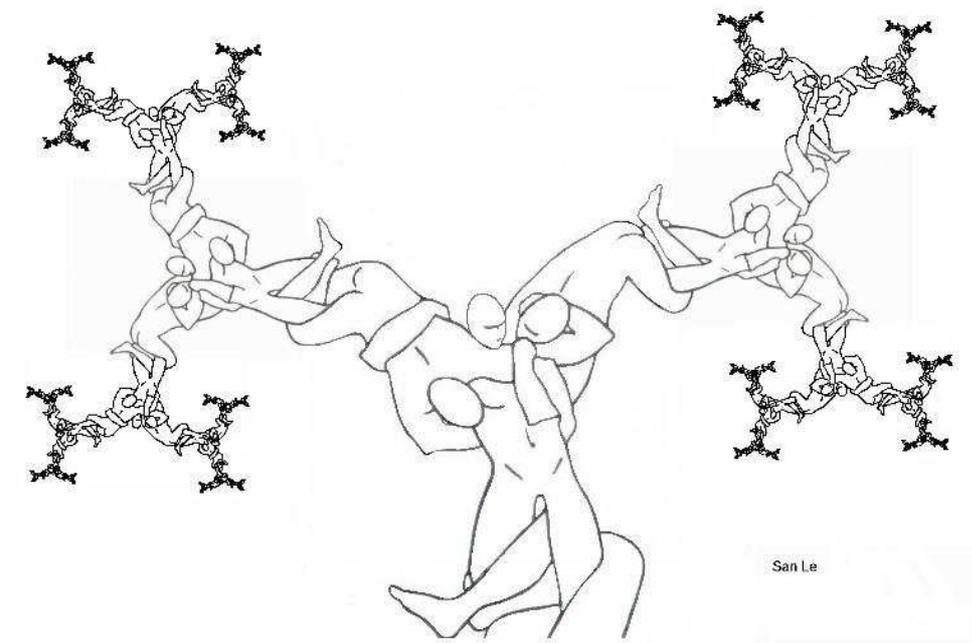}}
  \caption{ Final image using the tile of Figure \ref{cell2D_crop_border2} }
  \label{fractal_men_small}
  \end{figure}

For fractals, Escher's \emph{Circle Limit} series invokes the concept but it can also be
interpreted as a tessellation of the hyperbolic plane\cite[]{Adcock00}.
To make a fractal tiling that purely manifests the notion of self similarity, many of the previous techniques
are applied again with a few modifications.  Starting
with a single tile, a matching rule is created for how other tiles will be connected to its
sides.  An additional constraint is that the added tiles are decreased in size
by an amount that allows for growth of the fractal without the branches
intersecting.  Figure \ref{cell2D_crop_border2}
shows one rectangular tile where subsequent additions decrease in size by 1/2.  So
the next tile is half of the original, and its side \emph{A} will connect with \emph{A'}
and \emph{A''}.  To improve the esthetics, the two connecting sides are made to alternate
such that side \emph{A''} connects to the next tile after it is rotated 90 degrees clockwise
while side \emph{A'} connects to the mirror image of this tile and rotated 90 degrees
counter-clockwise.  Figure \ref{fractal_men_small} shows the result.

\begin{figure}
  \centerline{\includegraphics[height=3.752in,width=4in,angle=0]{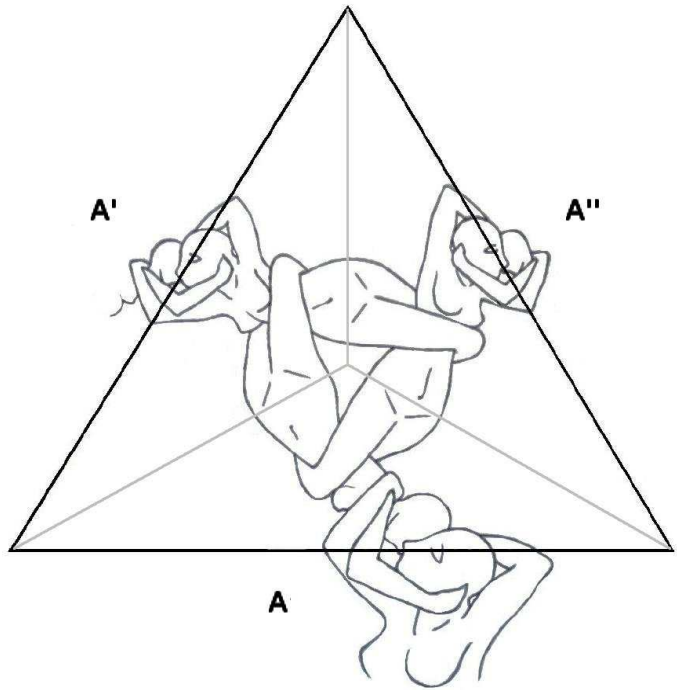}}
  \caption{ Fractal tessellation combination tile }
  \label{cellA_border_crop3}
  \end{figure}

\begin{figure}
  \centerline{\includegraphics[height=4.22in,width=5in,angle=0]{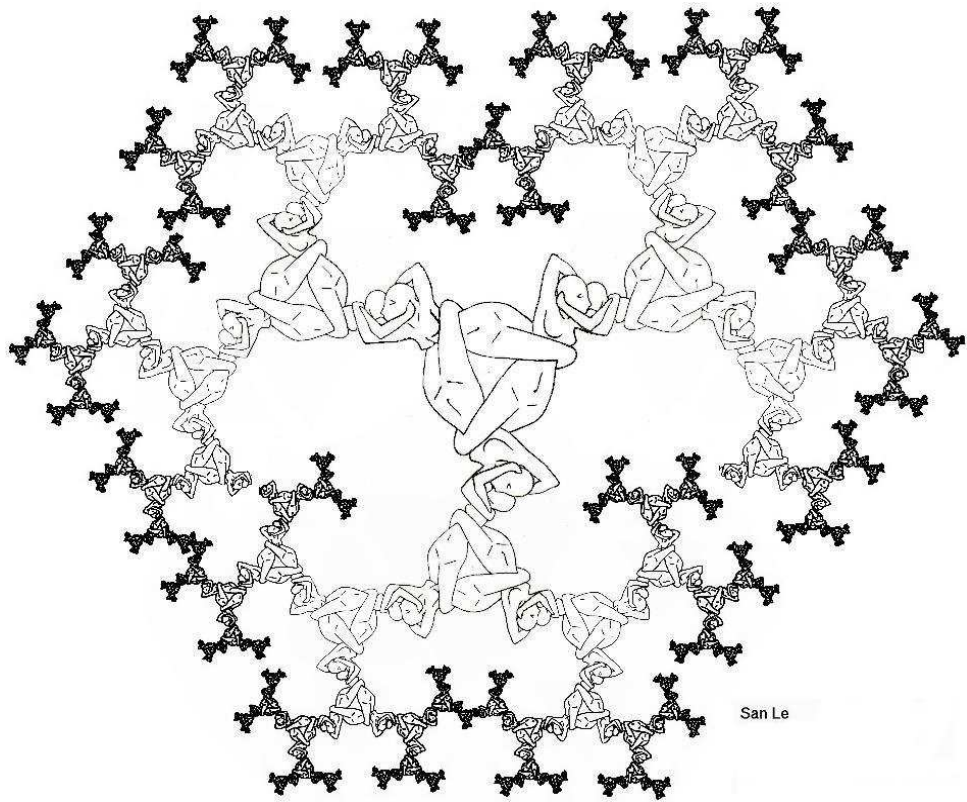}}
  \caption{ Final image using the tile of Figure \ref{cellA_border_crop3} }
  \label{fractal_women_small}
  \end{figure}

For the next fractal tile made of an isosceles triangle, the same mirroring connection on
side \emph{A'} is made.  In addition, the principals of tessellation are
used so that the fractal can progress in other directions.  As seen in \ref{cellA_border_crop3},
the area around the centre of the triangle is a region of symmetry for the design which is repeated every 
120 degrees.  This is delineated by 3 sub-triangles which form the tessellation inside the
isosceles.  From here, the lower triangle can be swapped with either of the other two leading to the
fractal of Figure \ref{fractal_women_small} which has branches in the downward direction.
In this fractal, the decrease of each tile was only one third rather than one half
which illustrates the problem of branches colliding.  Despite this, the final result
is very appealing.

\section{Concluding Remarks}

Escher's work introduced the world to the beauty of geometrical art.  But
non-mathematician artists tended not to follow his example, and so a wealth of
trigonometric shapes only exists as blank tiles waiting to be filled.
By describing the process of incorporating tessellations and fractals into
art, we hope to show that the challenges are artistic rather than mathematical.

In the same way that mathematical analysis can give a deeper appreciation of
art, Escher used art as a means of gaining insight into the mathematics of geometry.
Continuing this tradition, our paper demonstrates that the connections contained
within tessellations and fractals are best seen when combined with tiling art. 
With the limitless possibilities as to what can be put inside a tile, artists
are well suited to find the undiscovered pattern contained in each.

\label{lastpage}

\end{document}